\documentstyle[11pt]{article}
\begin{document}

\hfill{arch-ive/9602003} \\
\hfill{January, 1996} \\

\begin{center}

{\Large \bf Contractions of algebraical structures and \\
different couplings of Cayley--Klein and Hopf structures.} \\
\vspace{5mm}
{\bf N.A.Gromov}

\end{center}

\begin{center}
Department of Mathematics, Komi Research Center,
Russian Academy of Sciences, \\
 16700 Syktyvkar, Russia,
e-mail: gromov@omkomi.intec.ru.
\end{center}
\vspace{5mm}

\section{ Contractions of algebraical structures}

 A {\bf Lie algebra} $ L $ is a vector space $ V $ equipped
with a Lie bracket $ [ \ , \ ]: L \otimes L \rightarrow L, $
which is linear in both arguments, antisymmetric and satisfy the
Jacobi identity. So the structure of Lie algebra is as follows:
$ L=(V,[ \ , \ ])$. It is well known fact that the non--degenerate
linear transformations $ \Phi: V \rightarrow V, \; det{\Phi} \not = 0 $
of vector space $ V $ give in result the Lie algebra $ L^{'} $
which is  isomorphic to the initial one's $ L $. How to receive a new Lie
algebra from $ L $? The answer was given many years ago by E.In{\"o}n{\"u} and
E.P.Wigner \cite{IW}: it is necessary to use a degenerate transformation
of $ V $. Let us regard $ \Phi_{\epsilon}: V \rightarrow V $ such that
$ det \Phi_{\epsilon} \not =0 $ for $ \epsilon \not =0 $, but
$ det \Phi_{\epsilon}=0 $ for $ \epsilon=0 $. Then if the limit

\begin{center}
\vspace{1.5cm} \unitlength=0.5mm \special{em:linewidth 0.4pt}
\linethickness{0.4pt}

\begin{picture}(120.00,40.00)

\put(10.00,10.00){\makebox(0,0)[cc]{$ L $}}
\put(100.00,60.00){\makebox(0,0)[cc]{ $ L \otimes L $}}

\put(85.00,10.00){\vector(-1,0){70.00}}
\put(100.00,55.00){\vector(0,-1){37.00}}
\put(85.00,60.00){\vector(-1,0){70.00}}

\put(10.00,60.00){\makebox(0,0)[cc]{$ L $}}
\put(100.00,10.00){\makebox(0,0)[cc]{$ L^{'} \otimes L^{'} $}}
\put(10.00,55.00){\vector(0,-1){37.00}}

\put(60.00,17.00){\makebox(0,0)[cc]{ $ [ \ , \ ]^{'} $ } }
\put(60.00,67.00){\makebox(0,0)[cc]{ $ [ \ , \ ] $ }}
\put(5.00,35.00){\makebox(0,0)[r]{ $ \Phi_{\epsilon} $ }}
\put(105.00,35.00){\makebox(0,0)[l]{ $ \Phi_{\epsilon} \otimes
 \Phi_{\epsilon}$ }}

\end{picture}
\end{center}

\begin{equation}
[ \ , \ ]^{'}=\lim_{\epsilon \to 0}\Phi_{\epsilon}\circ [ \ , \ ]\circ
 (\Phi_{\epsilon}^{-1} \otimes \Phi_{\epsilon}^{-1})
\label{1}
\end{equation}
exist, we have the new Lie algebra $ L^{'}=(V,[ \ , \ ]^{'}) $
with the same underlying vector space. The limit (\ref{1}) is
called a {\it contraction of Lie algebra}.

{\bf Graded contractions} was introduced in \cite{P}, \cite{MP} and
deal with the structure $ (L,\Gamma) $, where $ L $ is a Lie
algebra and $ \Gamma $ is a finite Abelian group. A grading
of $ L $ by $ \Gamma $ means the decomposition of the vector space
$ V $ of $ L $ into a direct sum of subspaces

\begin{equation}
V=\bigoplus_{\mu \in \Gamma} V_{\mu}
\label{2}
\end{equation}
such that for every choice of $  X \in V_{\mu}, \  Y \in V_{\nu} $,
if $ [X,Y] \not =0 $, then $ [X,Y] \in V_{\mu+\nu} $.
Symbolically

\begin{equation}
0 \not = [V_{\mu},V_{\nu}] \subseteq V_{\mu+ \nu}, \  \mu, \nu, \mu+\nu \in
 \Gamma.
\label{3}
\end{equation}
In this case the map $ \Phi_{\epsilon} $ additionally to
existence of (\ref{1}) have to preserve the grading equations
(\ref{2}), (\ref{3}). Graded contractions leads to a structure
$ (L^{'}, \Gamma)=(V,[ \ , \ ]^{'}, \Gamma) $ with the same
underlying $ V $ and the same grading group $ \Gamma $.

 A {\bf Lie bialgebra} $ (L, \eta) $ is a Lie algebra $ L $ endowed with
a cocomutator  $ \eta:L \rightarrow L \otimes L $ such that

i) $ \eta $ is 1--cocycle, i.e.
\begin{equation}
\eta([X,Y])=[\eta(X), 1 \otimes Y+Y \otimes 1]+
[1 \otimes X+X \otimes 1, \eta(Y)], \ \forall X,Y \in L,
\label{4}
\end{equation}

ii) the dual map $ \eta^{*}:L^{*} \otimes L^{*} \rightarrow L^{*} $
is a Lie bracket on $ L^{*} $.

Contracted Lie brackets are obtained by (\ref{1}).
Contraction of cocommutator $ \eta $ may be realized by the mapping
$ \Psi_{\epsilon}:L \rightarrow L^{'} $, which is in general
different from $ \Phi_{\epsilon} $, as follows:

\begin{center}
\vspace{1cm} \unitlength=0.5mm \special{em:linewidth 0.4pt}
\linethickness{0.4pt}

\begin{picture}(140.00,40.00)

\put(10.00,16.00){\makebox(0,0)[cc]{$ L^{'} \otimes L^{'} $}}
\put(65.00,16.00){\makebox(0,0)[cc]{ $ L^{'} $}}
\put(125.00,16.00){\makebox(0,0)[cc]{$ L^{'} \otimes L^{'} $}}
\put(10.00,60.00){\makebox(0,0)[cc]{$ L \otimes L $}}
\put(65.00,60.00){\makebox(0,0)[cc]{ $ L $}}
\put(125.00,60.00){\makebox(0,0)[cc]{ $ L \otimes L $}}

\put(16.00,55.00){\vector(0,-1){35.00}}
\put(65.00,55.00){\vector(0,-1){35.00}}
\put(115.00,55.00){\vector(0,-1){35.00}}
\put(60.00,16.00){\vector(-1,0){36.00}}
\put(110.00,16.00){\vector(-1,0){36.00}}
\put(110.00,60.00){\vector(-1,0){36.00}}
\put(60.00,60.00){\vector(-1,0){36.00}}


\put(35.00,22.00){\makebox(0,0)[cc]{ $ {\eta}^{'} $ } }
\put(35.00,65.00){\makebox(0,0)[cc]{ $ \eta $ }}
\put(12.00,40.00){\makebox(0,0)[r]{ $ \Psi_{\epsilon} \otimes \Psi_{\epsilon}
 $ }}
\put(118.00,40.00){\makebox(0,0)[l]{     $    \Phi_{\epsilon}
 \otimes \Phi_{\epsilon} $ }}
\put(63.00,40.00){\makebox(0,0)[r]{ $ \Psi_{\epsilon} $ }}
\put(67.00,40.00){\makebox(0,0)[l]{ $ \Phi_{\epsilon} $ }}
\put(90.00,22.00){\makebox(0,)[c]{ $ [ \ ,  \ ]^{'} $ }}
\put(90.00,65.00){\makebox(0,)[c]{ $ [ \ ,  \ ] $ }}
\end{picture}
\end{center}

\begin{equation}
\eta^{'}=\lim_{\epsilon \to 0}(\Psi_{\epsilon} \otimes \Psi_{\epsilon})\circ
 {\eta}\circ {\Psi}-{\epsilon}^{-1}.
\label{5}
\end{equation}

The consistency of $ [ \ , \ ] $ and $ \eta, $
which means the completely commutativity of diagram,
 gives the another
expression for cocommutator
\begin{equation}
\eta^{'}=\lim_{\epsilon \to 0}(\Psi_{\epsilon} \otimes \Psi_{\epsilon})\circ
 {\eta}\circ {\Phi}^{-1}_{\epsilon}.
\label{6}
\end{equation}
If all limits (\ref{1}), (\ref{5}), (\ref{6}) exist, then
$ (L^{'},{\eta}^{'}) $ is {\it a contracted Lie bialgebra} \cite{VG}.

An associative algebra $ A $ is said to be a {\bf Hopf algebra} if there
exist maps: coproduct $ \Delta: A \rightarrow A \otimes A, $
counit $ u: A \rightarrow {\rm C} $ and antipode
$ \gamma: A \rightarrow A $, such that $ \forall a \in A $
 the following axiomas are hold

\begin{eqnarray}
H_1:  &  (id  \otimes   \Delta)\Delta(a)=(   \Delta   \otimes
 id)\Delta(a), \nonumber \\
H_2: &(id \otimes u)\Delta(a)=(u \otimes id)\triangle(a), \\
H_3: & m((id \otimes \gamma)\Delta(a))=m((\gamma \otimes id)\Delta(a)),
 \nonumber
\label{7}
\end{eqnarray}
where $ m(a \otimes b)=ab $ is the usual multiplication.
So the Hopf algebra structure is $ (A, \Delta, u, \epsilon)$.
In general contractions of coproduct, counit and antipode
are realized by a different mapping of $ A $. But the map
$ {\rm X}: A \rightarrow A^{'} $, connected with coproduct
$ \Delta $, is the essential one (at least in case of quantum
deformations of an universal enveloping algebras).
 Therefore, if the limits
\begin{center}
\vspace{1cm} \unitlength=0.4mm \special{em:linewidth 0.4pt}
\linethickness{0.4pt}

\begin{picture}(120.00,40.00)

\put(10.00,10.00){\makebox(0,0)[cc]{$ A^{'} $}}
\put(105.00,60.00){\makebox(0,0)[cc]{ $ A \otimes A $}}

\put(18.00,10.00){\vector(1,0){68.00}}
\put(100.00,55.00){\vector(0,-1){37.00}}
\put(18.00,60.00){\vector(1,0){68.00}}

\put(10.00,60.00){\makebox(0,0)[cc]{$ A $}}
\put(105.00,10.00){\makebox(0,0)[cc]{$ A^{'} \otimes A^{'} $}}
\put(10.00,55.00){\vector(0,-1){37.00}}

\put(60.00,18.00){\makebox(0,0)[cc]{ $ {\Delta}^{'} $ } }
\put(60.00,66.00){\makebox(0,0)[cc]{ $ \Delta $ }}
\put(5.00,35.00){\makebox(0,0)[r]{ $ \chi_{\epsilon} $ }}
\put(105.00,35.00){\makebox(0,0)[l]{ $ \chi_{\epsilon} \otimes
 \chi_{\epsilon}$ }}

\end{picture}
\end{center}

$$
\Delta^{'}=\lim_{\epsilon \to 0}({\chi}_{\epsilon} \otimes
 {\chi}_{\epsilon})\circ {\Delta}\circ {\chi}_{\epsilon}^{-1},
$$

\begin{center}
\vspace{1cm} \unitlength=0.4mm \special{em:linewidth 0.4pt}
\linethickness{0.4pt}
\begin{picture}(120.00,40.00)

\put(10.00,10.00){\makebox(0,0)[cc]{$ A^{'} $}}
\put(100.00,60.00){\makebox(0,0)[cc]{ $ {\rm C} $}}

\put(20.00,10.00){\vector(2,1){68.00}}
\put(18.00,60.00){\vector(1,0){68.00}}

\put(10.00,60.00){\makebox(0,0)[cc]{$ A $}}
\put(10.00,55.00){\vector(0,-1){37.00}}

\put(50.00,20.00){\makebox(2,1)[cc]{ $ {u}^{'} $ } }
\put(50.00,66.00){\makebox(0,0)[cc]{ $ u $ }}
\put(5.00,35.00){\makebox(0,0)[r]{ $ \chi_{\epsilon} $ }}

\end{picture}
\end{center}

\begin{equation}
u^{'}=\lim_{\epsilon \to 0}u\circ {\chi}_{\epsilon}^{-1},
\label{8}
\end{equation}

\begin{center}
\vspace{1cm} \unitlength=0.4mm \special{em:linewidth 0.4pt}
\linethickness{0.4pt}
\begin{picture}(120.00,40.00)

\put(10.00,10.00){\makebox(0,0)[cc]{$ A^{'} $}}
\put(100.00,60.00){\makebox(0,0)[cc]{ $ A $}}

\put(18.00,10.00){\vector(1,0){68.00}}
\put(100.00,55.00){\vector(0,-1){37.00}}
\put(18.00,60.00){\vector(1,0){68.00}}

\put(10.00,60.00){\makebox(0,0)[cc]{$ A $}}
\put(100.00,10.00){\makebox(0,0)[cc]{$ A^{'} $}}
\put(10.00,55.00){\vector(0,-1){37.00}}

\put(60.00,16.00){\makebox(0,0)[cc]{ $ {\gamma}^{'} $ } }
\put(60.00,66.00){\makebox(0,0)[cc]{ $ \gamma $ }}
\put(5.00,35.00){\makebox(0,0)[r]{ $ \chi_{\epsilon} $ }}
\put(105.00,35.00){\makebox(0,0)[l]{ $ \chi_{\epsilon} $ }}

\end{picture}
\end{center}
$$
\gamma^{'}=\lim_{\epsilon \to 0}{\chi}_{\epsilon}\circ {\gamma}\circ
 {\chi}_{\epsilon}^{-1}.
$$
exist, then $ (A^{'}, {\Delta}^{'}, u^{'}, {\epsilon}^{'}) $
is {\it a contracted Hopf algebra.}
Contractions of quantum groups (or Hopf algebras) was first developed in
\cite{C1}, \cite{C2}.

The notion of contraction first introduced by E.In{\"o}n{\"u}
and E.P.Wigner for Lie groups and algebras is now extend
on new algebraical structures. Some of them are mention in this
talk. But the fundamental idea of degenerate transformation
is presented in all cases.

\section{ Different couplings of Cayley--Klein and Hopf structures}

The universal enveloping algebra $ U(L) $ of Lie algebra $ L $ may be
defined as a Hopf algebra $ U_z(L) $, if $ \Delta, u, \gamma $ are introduced
 in
appropriate manner. For $ U_z(L) $ it is sufficient to define
 $ \Delta, u, \gamma $ only for generators $ \{ {\rm X}_k \} $ of $ L $.

The quantum algebra $ U_{z^*}(so(3)) \equiv so_{z^*}(3;X^*_{02}) $
in the rotation basis $ \{ X_{01},X_{02},X_{12} \} $, where
$ X_{\mu \nu}, \  \mu< \nu, \  \mu,\nu=0,1,2 $
are the generators of rotations in 2--planes
$ \{x_{\mu},x_{\nu} \},  \   x_{\mu} $ --- Cartesian coordinates
in 3--dimentional Euclidean space, have the following coproduct \cite{V1}
$$
\Delta X_{02}^*=1 \otimes X_{02}^{*}+X_{02}^{*} \otimes 1, \  X_{02}^{*}
\mbox{ -- primitive operator},
$$

\begin{equation}
\Delta X^{*}=X \otimes e^{{z^* \over 2}X_{02}^*}+
e^{-{z^* \over 2}X_{02}^*} \otimes X^{*}, \  X^*=X_{01}^*, X_{12}^*,
\label{9}
\end{equation}
counit $ \epsilon(X^*_{\mu, \nu})=0, $
antipode $ \gamma(X_{02}^*)=-X_{02}^* $,
\begin{equation}
\gamma(X_{01}^*)=-X_{01}^*cos{z^* \over 2}+X_{12}^*sin{\displaystyle{z^*
 \over 2}},
\label{10}
\end{equation}
$$
\gamma(X_{12}^*)=-X_{12}^*cos{z^* \over 2}-X_{01}^*sin{\displaystyle{z^*
 \over 2}},
$$
and deformed commutators
\begin{equation}
 [X_{01}^*,X_{02}^*]=X_{12}^*, \  [X_{02}^*,X_{12}^*]=X_{01}^*, \
[X_{12}^*,X_{01}^*]={\displaystyle{1 \over z^*}}sinh(z^*X_{02}^*)
\label{11}
\end{equation}

The orthogonal Cayley--Klein (CK) algebras $ so(n+1,
 {\bf j})=\{ X_{\mu,\nu} \} $
are obtained from $ so(n+1)=\{ X_{\mu\nu}^* \} $
by the tranformations \cite{G1}, \cite{G2}
\begin{equation}
X_{\mu\nu}=J_{\mu\nu}X_{\mu\nu}^*, \
J_{\mu\nu}=\prod^{\nu}_{k=\mu+1}{j_k}, \  \mu<\nu, \ \mu, \nu=0,1, \ldots, n,
\label{12}
\end{equation}
where parameters $ j_k=1,\iota_k,i, \  k=1,2, \ldots, n, $
and $ \iota_k $ are the dual (or Study) units with the algebraic
properties:
$ \iota_k^2=0, \  \iota_k \iota_r=\iota_r \iota_k \not =0, \ k \not =r, \
 \iota_k / \iota_k=1. $
The dual values of $ j_k $ correspond to In{\"o}n{\"u}--Wigner
 contractions \cite{IW}
and imaginary values correspond to analytical continuations
(or Weyl unitary trick).

The quantum CK algebras $ so_{z}(n+1;{\bf j}) $ are obtained
from quantum algebra $ so_{z^*}(n+1) $ by the transformations (\ref{12}),
supplemented with the transformation
\begin{equation}
z^*=Jz
\label{13}
\end{equation}
of the deformation parameter, where the multiplier $ J $
depends on the choice of the set of primitive generators
of quantum algebra.

Under introducing of a Hopf algebra structure in the universal
enveloping algebra some commuting generators of quantum
algebra $ so_{z^*}(n+1) $ have to be selected as the primitive
elements of the Hopf algebra. Let us observe, that any permutation
$ \sigma \in S(n+1) $ indices of the rotation generators
$ X^*_{\mu \nu} $ of $ so_{z^*}(n+1) $ (i.e. the transformations
from the Weyl group) leads to the isomorphic Hopf algebra
$  so^{'}_{z^*}(n+1) $, but with some other set of primitive
generators. This isomorphism of the Hopf algebras may be destroyed
by contractions. Indeed, under transformations (\ref{12}), (\ref{13})
the primitive generators of $ so_{z^*}(n+1) $ and
$  so^{'}_{z^*}(n+1) $ are multiplied by the different
$ J_{\mu \nu} $ and for the specific dual values of $ j_k $,
this may lead to the nonisomorphic quantum CK algebras.
In other words, Hopf and CK structures may be combined
in a different manner for quantum CK algebras.

Let us illustrate the different couplings of Hopf and CK structures
for the case of $ so_z(3;{\bf j}), \ {\bf j}=(j_1,j_2). $
The application of (\ref{12}) and (\ref{13}), which reads as
$ z^*=j_1j_2z, $ to (\ref{9})--(\ref{11}) gives the quantum CK
algebra $ so_z(3;{\bf j};X_{02}) $ with coproduct (\ref{9}), where
generators and deformation parameter are now without star,
counite: $ \epsilon(X_{\mu \nu})=0 $,
antipode: $ \gamma(X_{02})=-X_{02}, $
$$
\gamma(X_{01})=-X_{01}cos{j_1j_2}{z \over 2}+
 X_{12}{j_1j_2^{-1}}sin{j_1j_2}{\displaystyle{z \over 2}},
$$
\begin{equation}
\gamma(X_{12})=-X_{12}cos{j_1j_2}{z \over 2}-
 X_{01}j_1^{-1}j_2sin{j_1j_2}{\displaystyle{z \over 2}},
\label{14}
\end{equation}
and deformed commutators
\begin{equation}
 [X_{01},X_{02}]=j_1^2X_{12}, \  [X_{02},X_{12}]=j_2^2X_{01}, \
[X_{12},X_{01}]={\displaystyle{1 \over z}}sinh(zX_{02}).
\label{15}
\end{equation}

The substitution $ \sigma(0)=0, \  \sigma(1)=2, \  \sigma(2)=1 $
indices of generators in (\ref{9})--(\ref{11}) and transformations
(\ref{12}), (\ref{13}), which reads as $ z^*=j_1z $, gives the
quantum CK algebra $ so_z(3;j;X_{01}) $ with antipode:
$ \gamma(X_{01})=-X_{01},$
$$
\gamma(X_{02})=-X_{02}cos{j_1}{z \over 2}-
 X_{12}{j_1}sin{j_1}{\displaystyle{z \over 2}},
$$
\begin{equation}
\gamma(X_{12})=-X_{12}cos{j_1}{z \over 2}+
 X_{02}j_1^{-1}sin{j_1}{\displaystyle{z \over 2}},
\label{16}
\end{equation}
and deformed commutators
\begin{equation}
 [X_{12},X_{01}]=X_{02}, \  [X_{01},X_{02}]=j_1^2X_{12}, \
[X_{02},X_{12}]=j_2^2{\displaystyle{1 \over z}}sinh(zX_{01}).
\label{17}
\end{equation}
The substitution  $ \sigma(0)=1, \  \sigma(1)=0, \  \sigma(2)=2 $
in (\ref{9})--(\ref{11}) and introducing CK parameters by
(\ref{12}), (\ref{13}), which reads as $ z^*=j_2z $, leads to the
quantum CK algebra $ so_z(3;{\bf j};X_{12}) $ with antipode:
$ \gamma(X_{12})=-X_{12},$
$$
\gamma(X_{01})=-X_{01}cos{j_2}{z \over 2}-
 X_{02}{j_2^{-1}}sin{j_2}{\displaystyle{z \over 2}},
$$
\begin{equation}
\gamma(X_{02})=-X_{02}cos{j_2}{z \over 2}+
 X_{01}j_2sin{j_2}{\displaystyle{z \over 2}},
\label{18}
\end{equation}
and deformed commutators
\begin{equation}
 [X_{02},X_{12}]=j_2^2X_{01}, \  [X_{12},X_{01}]=X_{02}, \
[X_{01},X_{02}]=j_1^2{\displaystyle{1 \over z}}sinh(zX_{12}).
\label{19}
\end{equation}

This three quantum CK algebras are contracted to quantum Galilei
algebras by putting $ j_1=\iota_1, \  j_2=\iota_2 $
in Eqs.(\ref{14})--(\ref{19}), which gives:
$$
  so_z(3;\iota_1,\iota_2;X_{02}) : \quad
\gamma(X)=-X, \ X=X_{02},X_{01},X_{12},
$$
\begin{equation}
 [X_{01},X_{02}]=0, \  [X_{02},X_{12}]=0, \
[X_{12},X_{01}]={\displaystyle{1 \over z}}sinh(zX_{02});
\label{20}
\end{equation}
$$
 so_z(3;\iota_1,\iota_2;X_{01}) : \quad
\gamma(X_{01})=-X_{01}, \ \gamma(X_{02})=-X_{02},
$$
$$
\gamma(X_{12})=-X_{12}+{\displaystyle{z \over 2}}X_{02},
$$
\begin{equation}
 [X_{02},X_{12}]=0, \  [X_{01},X_{02}]=0, \
[X_{12},X_{01}]=X_{02}.
\label{21}
\end{equation}
$$
 so_z(3;\iota_1,\iota_2;X_{12}) : \quad
\gamma(X_{12})=-X_{12}, \ \gamma(X_{02})=-X_{02},
$$
$$
\gamma(X_{01})=-X_{01}-{\displaystyle{z \over 2}}X_{02},
$$
\begin{equation}
 [X_{02},X_{12}]=0, \  [X_{01},X_{02}]=0, \
[X_{12},X_{01}]=X_{02}.
\label{22}
\end{equation}

Mathematically quantum algebras (\ref{21}) and (\ref{22})
are isomorphic as Hopf algebras and both are not isomorphic
with quantum algebra (\ref{20}). But quantum algebras (\ref{21})
and (\ref{22}) are different from physical (kinematical) point
of view since the primitive operator $ X_{01} $ of (\ref{21})
is interpreted as generator of time translation while the
primitive operator  $ X_{12} $ of (\ref{22}) is interpreted
as generator of Galilei boost. As far as primitive operators
  are appear in coproduct (cf.(\ref{9}))
of (\ref{21}) and (\ref{22}) in the form
$ exp({\displaystyle{z \over 2}}X_{01}) $ and
$ exp({\displaystyle{z \over 2}}X_{12}), $
respectively, then the deformation parameters have the different
physical dimensions: $ [z]=sec $ for Galilei quantum algebra
(\ref{21}) and $ [z]=cm/sec $ for Galilei quantum algebra
(\ref{22}).

In the case of quantum algebra $ so_z(4;{\bf j}) $
there are 6 possibilities to choose two primitive generators
$ X_{\mu\nu} $ and $ X_{{\mu}'{\nu}'}, \
 \mu \not = \nu \not = {\mu}' \not = {\nu}'. $
It turn out that for some choice not all CK contractions are
allowed. But there is such pair of primitive operators
$ X_{03}, X_{12}, $ when all CK contractions are allowed.
In the next case of $ so_z(5;{\bf j}) $ there are 30
different choice of two primitive operators.
For some choice all CK contractions are forbidden.
All restrictions on contractions are arisen from coproduct.
But again for primitive generators $ X_{04}, X_{13} $
all CK contractions are allowed. Moreover the following
general statement holds: for any quantum algebra
$ so_z(n+1;{\bf j}), \  {\bf j}=(j_1, \ldots, j_n) $
there exist the set of primitive operators
$ X_{0n}, X_{1,n-1}, \ldots ,X_{k-1,k}, \  k=[(n+1)/2] $,
when all CK contractions $ j_m=\iota_m, \  m=1, \ldots,n $
are allowed.

So suggested constructive algorithm may be
applied for description physical systems with
symmetries of nonsemisimple quantum algebras
by the help of contraction physical systems
(usually well known) with symmetries of
simple quantum algebras in such a way that a set
of primitive operators of resulting system may have
different physical interpretations.

\end{document}